\documentstyle[openbib,twoside,titlepage]{article}

\newcommand{\ttl}{Recent Contributions to the Calculus of Finite
Differences: A Survey}

\newcommand{\aug}[1]{\left\langle #1 \right\rangle}
\newcommand{\r}[1]{\left\lfloor #1 \right\rceil}
\newcommand{\rc}[2]{\fee \r{\begin{array}{c} #1 \\ #2 \end{array}} \foo}
\newcommand{\D}{{\mbox{D}}}
\newcommand{\infnty}{\infty}

\newtheorem{prop}{Proposition}
\newtheorem{thm}{Theorem}
\newtheorem{cor}{Corollary}
\newtheorem{lem}{Lemma}[thm]
\newcommand{\proof}{{\em Proof:} }

\newcommand{\foo}{\renewcommand{\baselinestretch}{1}
        \setlength{\parskip}{10pt plus 5pt minus 5pt}}
\newcommand{\fee}{\renewcommand{\baselinestretch}{1}
        \setlength{\parskip}{0pt}}
\foo

\begin{document}
\title{\ttl}
\author{D Loeb\thanks{Supported by a Chateaubriand Fellowship}\\
8 Rue Claude Terrasse\\ 75016 Paris FRANCE
\and G-C Rota\thanks{Supported by National Science Foundation under
Grant MCS-8104855}\\
MIT Dept. of Mathematics\\ Cambridge, MA 02139 USA}
\date{Copyright \today}
\maketitle

\begin{abstract}
We retrace the recent history of the Umbral Calculus. After studying
the classic results concerning polynomial sequences of binomial type,
we generalize to a certain type of logarithmic series. Finally, we
demonstrate numerous typical examples of our theory. 

{\bf  Une Pr\'{e}sentation des Contributions Recents au Calcul des
Dif\-f\'{e}\-rences Finies}

Nous passons en revue ici les resultats recents du calcul ombral. Nous
nous interessons tout d'abord aux resultats classique appliqu\'{e}s
aux suites de polyn\^omes de type binomial, pius elargions le champ
d'\'{e}tude aux series logarithmiques. Enfin nous donnons de nombreaux
exemples types d'application de cette th\'{e}orie.
\end{abstract}
\tableofcontents

\section{Polynomials of Binomial Type}
Much of the calculus of finite differences leans upon certain analogies
between two linear operators defined on functions $f(x)$ ranging over
a field $K$ of characteristic zero, and taking values in $K$:
\begin{enumerate}
\item The {\em derivative operator} \D, defined as
$$ \D f(x) = f'(x), $$
and
\item The {\em forward difference operator} $\Delta$, defined as
$$ \Delta f(x) = f(x+1) - f(x). $$
\end{enumerate}

\subsection{The  Powers of $x$}
The sequence of polynomials\footnote{When one speaks of a sequence of
polynomials, one generally means a sequence $p_n(x)$ indexed by a
nonnegative
integer $n$ such that the degree of each polynomial $p_n(x)$ is
exactly equal to $n$.} $p_n(x)=x^n$ is related
to the derivative operator by the elementary property
\begin{equation}\label{star}
\D p_n(x) = n p_{n-1}(x).
\end{equation}
Such sequences of polynomials have been called {\em Appell}. The
sequence $x^n$ is
distinguished  among all Appell sequences by the further property
$$ p_n(0) = \delta_{n0}$$
for all $n\geq 0$.

Furthermore, $p_n(x)=x^n$ plays a key role in Taylor's theorem
$$f(x+a) = \sum_{n=0}^{\infty} p_n(a)\frac{\D ^n f(x)}{n!}
 = \sum_{n=0}^{\infty} \frac{a^n \D ^n f(x)}{n!}$$
for all $a\in K$ and for all suitable functions $f(x)$.

Finally, we note that the sequence $p_n(x)$ is of {\em binomial type}.
That is to say, it satisifies the {\em binomial theorem}
\begin{equation}\label{binomial}
p_n (x+a) = \sum_{k=0} ^{\infty} {n \choose k } p_k(a) p_{n-k}(x)
\end{equation}
for all $a\in K$ and all $n\geq 0$.

\subsection{The  Lower Factorial}

Define another sequence of polynomials called the {\em lower factorial
 powers} as follows
$$p_0(x)=(x)_0=1$$
and
$$p_n(x)=(x)_n = x(x-1)(x-2)\cdots(x-n+1)$$
for $n\geq 0 $. It can be interpretted as the number of {\em
injections} from an $n$-element set to an $x$-element set.\cite{sym}

The lower factorial is related to the forward difference operator  by
the property
$$ \Delta p_n(x) =  np_{n-1}(x).$$
Such sequences are called {\em Sheffer sequences} for the forward difference
operator. The lower factorial is distinguished among all Sheffer
sequences by the further property
$$ p_n(0)= \delta_{n0}$$
for $n\geq 0$.

Furthermore, for suitable functions $f(x)$ and all constants $a\in K$, we
have {\em Newton's formula}
$$f(x+a) = \sum_{n=0}^{\infty} (a)_n \frac{\Delta^n f(x)}{n!}.$$

Lastly,  {\em Vandermonde's identity} states that the sequence
$p_n(x)=(x)_n$ is of binomial type:
$$(x+a)_n = \sum_{k=0}^\infty {n \choose k} (a)_k
(x)_{n-k}. $$

\subsection{The General Case}
Our initial problem will be that of carrying this analogy as far as
reasonably possible in order to make the parallel of $\D$
and $\Delta$ into an accidental special case.

To this end, we begin by classifying all sequences of polynomials
of binomial type, that is, all sequences of polynomials $p_n(x)$
(with $\deg(p_n(x))=n$ satisfying eq.~\ref{binomial}
$$p_n (x+a) = \sum_{k=0} ^{\infty} {n \choose k } p_k(a) p_{n-k}(x).$$

\begin{prop}\label{zero} A sequence of polynomials $p_n(x)$ of binomial type
satisfies
$$p_n(0)= \delta_{n0}$$
for all $n\geq 0$.
\end{prop}

\proof{By eq.~\ref{binomial} with $n=a=0$, $p_0(0)=1$. Now, for $n$
positive,
$$ 0=p_n(x)-p_n(x)=\sum_{k=1}^{n}{n \choose k} p_k(0) p_{n-k}(x).$$
However, $p_{k}(x)$ is a basis for the vector space of polynomials,
and in characteristic zero the binomial coefficients indicated are
never zero, so $p_k(0)=0.\Box$}

Now, to each sequence of binomial type $p_n(x)$, we associate a linear
operator $Q$ defined as
\begin{equation}\label{basic}
Qp_n(x) = np_{n-1}(x).
\end{equation}
For the moment, the operator $Q$ is defined only for polynomials;
later (Section~\ref{general}), we shall extend its domain to all more
general types of functions.

By iteration, we obtain
$$ Q^k p_n(x) = (n)_k p_{n-k} (x).$$
Recalling that ${n \choose k} = (n)_k/k!$, we can write the binomial
identity (eq.~\ref{binomial}) in the form
$$ p_j (x+a) = \sum_{k=0}^{\infnty} p_k(a) \frac{Q^k p_j(x)}{k!}$$
for $j\geq 0$.

Now, let $p(x)$ be any polynomial. Given its expansion
$p(x)=\sum_{j=0}^\infty c_j p_j(x) $ in terms of the basis $p_j(x)$,
we  have
\begin{eqnarray}
p(x+a) & = & \sum_{k=0}^\infnty c_j p_j(x+a) \nonumber \\
&=& \sum_{k=0}^\infnty p_k(a) Q^k \left( \sum_{j=0}^\infnty
c_j p_j(x) \right) \nonumber \\
&=& \sum_{k=0}^\infnty p_k(a) \frac{Q^k}{k!} p(x). \label{import}
\end{eqnarray}
This calculation for all polynomials $p(x)$ can be recast into a more
elegant form by introducing the
{\em shift operator}
$$ E^a f(x) = f(x+a). $$
We then can write eq.~\ref{import} as the operator identity
\begin{equation}\label{taylor}
E^a = \sum_{k=0}^\infnty p_k(a) Q^k/k!
\end{equation}
which we shall call the {\em Taylor's Formula for the sequence
$p_n(x)$.}
Hence, Newton's formula is nothing more than Taylor's formula for the
lower factorials.
Again, eq.~\ref{taylor} has been proven only for operations on
polynomials for the time being; however, we shall see (Theorem
\ref{tailor}) that its domain of validity is considerably greater.

The right side of eq.~\ref{taylor} obviously commutes with $Q$. Thus,
$$QE^a = E^aQ $$
for all constants $a$. Any linear operator $T$ with the property that
\begin{equation}\label{eat}
TE^a = E^aT
\end{equation}
for all scalars $a$ will be said to be a {\em shift-invariant operator}. A
linear shift-invariant operator is called a {\em delta  operator} if its
kernel is exactly the field of all constants $K$.

We have thus proved all of Theorem~\ref{1a} and the hard part of
Theorem~\ref{1b}.
\begin{thm}\label{1a}
Let $p_n(x)$ be a sequence of polynomials of binomial type. Then
there exists a unique delta operator such that
$$Qp_n(x)= np_{n-1}(x)$$
for $n\geq 0.\Box$
\end{thm}

Next, we have the converse.
\begin{thm}\label{1b}
Let $Q$ be a delta operator. Then there exists a unique sequence of
polynomials of binomial type $p_n(x)$ satisfying
$$Qp_n(x)= np_{n-1}(x)$$
for $n\geq 0$.
\end{thm}

\proof{We can recursively define $p_n(x)$ subject to the requirements
that
\begin{itemize}
\item $p_0(x)=1$,
\item For $n$ positive, $p_n(x)$ is in the inverse image of $Q$ on
$np_{n-1}(x)$, and
\item For $n$ positive, $p_n(0)=0$.
\end{itemize}
It remains only to show that the resulting sequence of polynomials
obeys eq.~\ref{binomial}. By the reasoning before Theorem~\ref{1a}, we again
have eq.~\ref{import}. Finally, setting $p(x)=p_n(x)$ yields the
binomial theorem.$\Box$}

Thus, we see that all of the features of the analogy between the
operators $\D$ and $\Delta$, and the polynomial sequences $x^n$ and
$(x)_n$ are shared by all other delta operators, and their
corresponding sequences of binomial type.

Sequences of polynomials of binomial type are of frequent occurence in
combinatorics,\footnote{In combinatorics, sequences of binomial type
enumerate  the number of functions from an
$n$-element set to an $x$-element set enriched with some
mathematical structure on each fiber. For details see \cite{sym}.}
probability, statistics, function theory and representation theory.

\subsection{The Abel Polynomials}\label{able}
We note the following striking example of a sequence of polynomials of
binomial type.
Given a constant $b$, define {\em Abel's polynomials}\/\footnote{For
example,  the Abel polynomial $A_n(x;-1)$ counts
the number of {\em reluctant functions} from an $n$-element set to an
$x$-element set. A reluctant function is a function enriched with a
rooted labelled forest on each fiber.} as follows
$$A_0(x;b)=1$$
and for $n$ positive
$$A_n(x;b)= x(x-nb)^{n-1}.$$

 It is easily verified that, setting $Q=E^b\D =\D E^b$, we have
$QA_n(x;b) = nA_{n-1}(x;b).$
Hence, Theorem~\ref{1b} implies at once {\em Abel's identity}
$$(x+a)(x+a-nb)^{n-1} = \sum_{k=0}^\infty {n \choose k} a(a-kb)^{k-1}
x(x-(n-k)b)^{n-k-1}.$$

Furthermore, Taylor's formula for the sequence of Abel polynomials
gives
$$p(x+a)= \sum_{k=0}^\infty \frac{a(a-kb)^{k-1}}{k!}D^kp(x+kb).$$
For example, for $p(x)=x^n$, we obtain
$$ (x+a)^n = \sum_{k=0}^\infty \frac{a(a-kb)^{k-1}}{k!}(n)_k(x+kb)^{n-k}$$
which is not {\em a priori} obvious.

In order to devise more examples, it will be helpful to first study
shift-invariant operators in greater detail.

\section{Shift-Invariant Operators}

For convenience, we introduce the following bit of notation. For any
function $f(x)$, we write
\begin{equation}\label{augment}
 \aug{f(x)}_{(0)} = f(0) = [f(x)]_{x=0}.
\end{equation}
The linear
functional $\aug{}_{(0)}$ is called the {\em augmentation} or {\em
evaluation at zero.}

Let $p_n(x)$ be a sequence of binomial type and let $Q$ be its delta
operator.  Taylor's theorem for $p_n(x)$ can be written in terms
of the augmentation (interchanging the roles of $x$ and $a$)
$$E^a p(x) = \sum _{k=0}^\infty p_k(x) \frac{\aug{E^a Q^k p(x) }_{(0)}}
{k!}. $$
Now, let $T$ be any shift-invariant operator, and $p(x)$ any
polynomial. We have
\begin{eqnarray*}
E^a T p(x) &=& T E^a p(x) \\
&=& \sum_{k=0}^\infty T p_k(x) \aug{E^a Q^k p(x) }_{(0)}/ k! .
\end{eqnarray*}
Again exchanging the role of $x$ and $a$, therefore
$$Tp(x+a)= \sum_{k=0}^\infty \frac {\aug{E^a T p_k(x)}_{(0)}} {k!} Q^k p(x).
$$
Setting $a=0$, we obtain the operator identity
\begin{equation}\label{unique}
T = \sum_{k=0}^\infty c_k Q^k / k!
\end{equation}
where $c_k = \aug{ T p_k(x)}_{(0)}$.

Thus, we have proven
\begin{thm}[Expansion Theorem]\label{expansion}
Given a shift-invariant operator $T$ and a delta operator $Q$
associated to a sequence of binomial type $p_n(x)$ there is a unique
expansion given by eq.~\ref{unique}.$\Box$
\end{thm}
Actually, some mild continuity conditions must be imposed, but we will
refer the reader to the bibliography for  such technical (and
ultimately trivial) questions of topology.

The expansion theorem is to operators essentially what the generalized
Taylor's theorem is to polynomials. Together they allow the expansion
of either operators or polynomials any basis.

The above result can be restated in more elightening terms. Let
$K[[t]]$ be the ring of all formal series (sometimes called {\em
Hurwitz series})  in the variable $t$ with coefficients in $K$
$$f(t)=c_0 + c_1t + c_2t^2/2! + c_3t^3/3! + \cdots .$$
Then we have
\begin{thm}[Isomorphism Theorem]\label{ice}
\begin{enumerate}
\item For every delta operator $Q$, one has the following isomorphism of the
ring of formal power series $K[[t]]$ onto the ring of shift-invariant
operators
\begin{equation}\label{two2}
f(t)\mapsto \sum_{k=0}^\infty c_k Q^k /k! =f (Q).
\end{equation}
\item Furthermore, let $\phi$ be an isomorphism of the ring of formal
power series with the ring of shift-invariant operators. Then there
exists a delta operator $Q$ such that the isomorphism $\phi$ is of
the form eq.~\ref{two2}.
\end{enumerate}
\end{thm}

\proof{Part one is immediate from the Expansion Theorem. Part two is
classical result concerning formal series where $Q=\phi(t).\Box$}

In particular for $Q= \D$, we notice that all shift-invariant
operators are formal power series in the derivative. That is, the ring
of shift-invariant operators is $K[[\D ]]$. Moreover, delta operators
$T=T(\D )$ are formal power series in the derivative inwhich the
coefficient of $\D ^0$ is zero, and the coefficient of $\D ^1$ is
nonzero. Such series are known as {\em delta series}.

For example, the shift operator is given by the formal power series
$$E^a = \exp (a\D ) = \sum_{k=0}^\infty a^n \D ^n /n!. $$
As another example, let $K={\bf C}$ be the complex field, and define
the {\em Laguerre operator} by
\begin{equation}\label{seven}
Lf(x) = -\int_0^\infty e^{-y} f'(x+y) dy.
\end{equation}
The Lagerre operator is obviously a shift-invariant operator since
$$ LE^af(x)=-\int_0^\infty e^{-y} f'(x+a+y) dy=E^aLf(x).$$
Moreover,
\begin{eqnarray*}
\aug{Lx^n}_{(0)} &=& -\int_0^\infty e^{-y} ny^{n-1} dy\\
&=& \left[e^{-y}y^n\right]_{y=0}^\infty
-\int_0^\infty e^{-y} ny^{n-1} dy\\
&=& -\int_0^\infty e^{-y} ny^{n-1} dy\\
&=& -\int_0^\infty e^{-y} n(n-1)y^{n-2} dy\\
& \vdots & \vdots \\
&=& -\int_0^\infty e^{-y} n!dy\\
&=& -n!,
\end{eqnarray*}
and hence by the Expansion Theorem (Theorem~\ref{expansion}), we infer
that $L$ is a delta operator and that
\begin{eqnarray}
L &=& -\D - \D^2 - \D^3 - \cdots \nonumber \\
&=& \D / ( \D - I ) \label{september}
\end{eqnarray}
where $I=\D ^0$ is the identity operator. We will later (Section
\ref{soon}) compute the  sequence of binomial type associated with the
Laguerre operator.

As another example, let $W$ be the {\em Weierstrass operator} defined
as
$$ Wf(x) = \frac{1}{\sqrt{2\pi}} \int_{-\infty}^\infty f(y)
\exp(-(x-y)^2/2) dy. $$
Again, we calculate $c_n=\aug{ L x^n}_{(0)} $.
First, note that for $n$ odd, the integrand in question is an odd
function, so the
integral from $-\infty$ to $\infty$ is zero. For $n=0$, it is easier
to compute $c_0^2$.
\begin{eqnarray*}
c_0^2 &=& \frac{1}{2\pi} \int_{-\infty}^\infty \int_{-\infty}^\infty
\exp(-(x^2+y^2)/2) dxdy\\
 &=& \frac{1}{2\pi} \int_{0}^{2\pi} \int_{0}^\infty
\exp(-r^2/2) rdrd\theta\\
 &=& \int_{0}^\infty \exp(-r^2/2) rdr\\
 &=& \int_{0}^\infty \exp(-u) du\\
 &=& 1.
\end{eqnarray*}
In general, by integration by parts, we have the following recurrence
for $c_n$
\begin{eqnarray*}
c_n &=&  \frac{1}{\sqrt{2\pi}} \int_{-\infty}^\infty y^n
\exp(-y^2/2) dy \\
&=& \frac{1}{\sqrt{2\pi}} \left(
-\left[\exp(-y^2/2)y^{n-1}\right]_{-\infty}^\infty + (n-1)
\int_{-\infty}^\infty y^{n-2} \exp(-y^2/2) dy \right) \\
&=& 0 + (n-1)c_{n-2}.
\end{eqnarray*}
Thus, for $n$ even,
$$ c_n = (n-1)(n-3)\cdots 5 \cdot 3 \cdot 1 = \frac{n!}{(n/2)!2^n}.$$
so that
$$W= \sum_{n=0}^\infty \frac{\D^{2n}}{n!2^n} = \exp (\D ^2/2).$$

From  the Isomorphism Theorem (Theorem~\ref{ice}), it follows that a
shift-invariant operator has a unique inverse if and only if $T_1\neq
0$. For example, the Weierstrass operator $W$ has a unique inverse.

Lastly, we compute the {\em Bernoulli operator}
$$Jf(x)=\int_x^{x+1}f(y)dy$$
where we have
$$ \aug{ J x^n}_{(0)}  = \int_0^1 y^ndy =
\left[\frac{y^{n+1}}{n+1}\right]_{y=0}^1 = \frac{1}{n+1}.$$
Thus,
$$J = \sum_{k=0}^\infty \frac{\D^k}{(k+1)!} = \frac{e^{\D} - I }{\D} =
\frac{\Delta}{\D} .$$

We close this section with a result which may be considered to be
fundamental.
\begin{thm} {\bf (Fun\-da\-men\-tal Theo\-rem of the Cal\-cu\-lus of
Fi\-nite Diff\-er\-ences)}
Let $p_n(x)$ be the sequence of binomial type associated to a delta
operator $Q(\D)$. Let $Q^{(-1)}(t)$  be the inverse
formal power series of $Q(t)$. That is, suppose that
$Q(Q^{(-1)}(t))=Q^{(-1)}(Q(t)) = t$. Then the exponential generating
function for $p_n(x)$ is
$$\sum_{n=0}^\infty p_n(x)t^n/n! = \exp (xQ^{(-1)}(t)).$$
\end{thm}

Since this result is so important, we give two proofs. A third proof
found in \cite{sym} is more powerful and applies equally to symmetric
functions.

{\em Proof 1:} By Proposition~\ref{zero} and eq.~\ref{basic}, the
left hand side is characterized by the property
$$ \aug{ Q(\D) ^k p(x,t)}_{(0)} = t^k ,$$
so it will suffice to verify the right hand side obeys this property
equally. Let $q_{nk}$ denote the coefficients of
$Q(\D) ^k=\sum_{n=0}^\infty q_{nk}\D^n/n!.$
Then
\begin{eqnarray*}
\aug{ Q(\D) ^k \exp (xQ^{(-1)}(t))}_{(0)} &=&
\sum_{n=0}^\infty\aug{ Q(\D) ^k x^n}_{(0)} (Q^{(-1)}(t))^n/n!\\
&=& \sum_{n=0}^\infty q_{nk} (Q^{(-1)}(t))^n/n!\\
&=& Q(Q^{(-1)}(t))^k\\
&=& t^k.\Box
\end{eqnarray*}

{\em Proof 2:} By Taylor's Formula (eq.~\ref{taylor}),
$$\exp(a\D) = \sum_{k=0}^\infnty p_k(a) Q^k/k!.$$
Substituting $x$ for $a$, and $Q^{(-1)}(t)$ for $\D$, we have the
desired result.$\Box$

\section{Extension to Logarithmic Series}\label{general}
\subsection{The Harmonic Logarithms}
We shall now extend the domain of every shift-invariant operator to
a more general domain of formal series originating from
the Hardy field\footnote{No knowledge of Hardy fields nor of
asymptotic expansions is expected of the reader.}
 ${\cal L}$ called the {Logarithmic Algebra} consisting of all
expansions of
real functions in a neighborhood of infinity in terms of the monomials
$$\ell^\alpha = x^{\alpha_0} \log(x)^{\alpha_1}
\log(\log(x))^{\alpha_2} \cdots$$
for all vectors of integers\footnote{We could achieve far greater
generality by allowing $\alpha$ to be a vector of reals as in
\cite{ILA1}. The definitions in \cite{SGE} provide one with the basic
tools to carry out these calculations.} where only a
finite\footnote{Since all the vectors we deal with here have finite
support, we will usually adopt the convention of not writing the
infinite sequence of zeroes they all end with. Thus, we write (1) for
the vector (1,0,0,0,\ldots).} number of the $\alpha_i$ are different
from zero.

That is $\cal L$ is equivalent to the set of all formal sums
$$\sum_\alpha b_\alpha \ell^\alpha$$
where the sum is over vectors of integers $\alpha$ with finite
support, and for all nonnegative integers $n$, all integers
$\alpha_0, \alpha_1, \ldots, \alpha_{n-1}$, and all integers $\beta$,
there exists finitely many $\alpha_n > \beta$ such that there exists
integers $\alpha_{n+1},\alpha_{n+2},\ldots$ (only finitely many
different from zero) such that $b_\alpha \neq 0$.

Another characterization of $\cal L$ \cite[Theorem 5.4A.4]{thesis} is
that $\cal L$ contains a dense subset which is the smallest proper
{\em field} extensions of the ring of polynomials with real
coefficients such that the derivative $\D$ is a derivation of the
entire field, and the antiderivative $\D^{-1}$ is well defined up to a
constant of integration.

It will be noted that the derivative of such an expression is
awkward
$$\D \ell^{\alpha}=
\alpha_0 \ell ^{^{(\alpha_0 -1, \alpha_1, \alpha_2, \alpha_3 \ldots)}} +
\alpha_1 \ell ^{^{(\alpha_0 -1, \alpha_1 -1, \alpha_2, \alpha_3 \ldots)}}+
\alpha_2 \ell ^{^{(\alpha_0 -1, \alpha_1 -1, \alpha_2 -1, \alpha_3 \ldots)}} +
\cdots .$$
Moreover, the formula for the antiderivative of such a function is not
known in general,
and even in the  case $\alpha_0\neq 1$ when it is known, it is given
by a horrendous
expression \cite[p. 87]{thesis}. Since the monomials $\ell^\alpha$ are
so unwieldy and do not even begin to compare to the powers of $x$ in
terms of ease and simplicity of calculation, we are forced to resort to
another basis for the vector space ${\cal L}$ which as will be seen is
the true {\em
logarithmic analog} of the sequence $x^n$. This basis is called the {\em
harmonic logarithm} $\lambda_n^\alpha(x)$.

For $n$ an integer, and $\alpha=(\alpha_1, \alpha_2, \ldots )$ a
vector of integers with only a finite number of nonzero entries, the
{\em harmonic logarithm} of degree $n$ and order
$\alpha$ is defined by the series
$$\lambda_n^\alpha (x) = \sum_{\rho}
\frac{s(-n,\rho_1)}{ \r{-n}!} \left[ \prod_{i=1}^\infty e_{\rho_i -
\rho_{i+1}} (\alpha_i, \alpha_i -1, \ldots, \alpha_i - \rho_i +1)\right]
\ell ^{(n),(\alpha-\rho)}$$
where
\begin{enumerate}
\item The {\em linear partition} $\rho= (\rho_1, \rho_2,\ldots)$ is a
 weakly decreasing  sequences
of nonnegative integers,
\item The {\em Roman factorial}\/\footnote{After Steve Roman.
\cite{AFS,AFS2,GBC}.} $r{n}!$ is given by the formula
$$ \r{n}! = \left\{ \begin{array}{ll}
                    n! & \mbox{for $n$ a nonnegative integer, and}\\
                    (-1)^{n-1}/(-n-1)! & \mbox{for $n$ a negative integer,}
                    \end{array} \right. $$
\item For any integer $n$ negative, zero or positive, and for any
nonnegative integer $k$, the {\em Stirling number} $s(n,k)$ is the
coefficient of $y^k$ in
the Taylor series expansion of $(y)_n=\Gamma(y+1)/\Gamma(y-n+1)$, and
\item The {\em elementary symmetric function} $e_n(x_1,x_2,\ldots)$ is defined
the generating function
$$ \sum_{n=0}^\infty e_n(x_1,x_2,\ldots)y^n = \prod_{k=1}^\infty (1+x_ky).$$
\end{enumerate}

It may be helpful to point out several important special cases.
\begin{enumerate}
\item The harmonic logarithms of degree zero are given by
$$\lambda_0^\alpha = \ell^{(0),\alpha}=(\log x)^\alpha_1 (\log \log
x)^\alpha_2 \cdots.$$
\item The harmonic logarithm of order $\alpha=(-1)$ and degree one
$\lambda_1^{(-1)}(x)$ is the
{\em logarithmic integral} $\mbox{li}(x)=\int_0^x dt/\log t .$
\item The harmonic logarithms of order $\alpha=(0)$ are given by
$$\lambda_n^{(0)}(x)= \left\{ \begin{array}{ll}
        x^n & \mbox{for $n\geq 0$, and}\\
        0 & \mbox{for $n<0$.}
        \end{array} \right. $$
\item The harmonic logarithms of order $\alpha=(1)$ are given by
$$\lambda_n^{(0)}(x)= \left\{ \begin{array}{ll}
        x^n\left(\log   x - 1 - \frac{1}{2} -  \cdots -
        \frac{1}{n}\right) & \mbox{for $n\geq 0$, and}\\
        x^n & \mbox{for $n<0$.}
        \end{array} \right. $$
\item The harmonic logarithms of order $\alpha=(2)$ are given by
$$\lambda_n^{(0)}(x)= \left\{ \begin{array}{ll}
         x^n (\log x)^2 -
        x^n\left(\frac{2}{1}+\frac{2}{2}+\cdots+\frac{2}{n}\right)  \\[-1mm]
         +x^n\left[\frac{2}{2}\left(1+\frac{1}{2}\right) + \cdots +
        \right. & \mbox{for $n\geq 0$, and}\\[-1mm]
        \left.
\frac{2}{n}\left(1+\frac{1}{2}+\cdots+\frac{1}{n}\right)\right] \\
        2x^n \left[\log(x)-1-\frac{1}{2}-\cdots
        -\frac{1}{n-1}\right].& \mbox{for $n<0$.}
        \end{array} \right. $$
\item The harmonic logarithms of order $\alpha=(t)$ where $t$ is a nonnegative
integer\footnote{Our original paper \cite{LA} dealt exclusively with
this case.}  are given by the following formulas
\begin{eqnarray*}
\sum_{t=0}^\infty \lambda_n^{(t)}(x)z^t/t! &=& \r{n}! x^{n+z}
\Gamma(z+1) / \Gamma(z+n+1)\\
\lambda_n^{(t)} &=& x^n \r{n}! (x\D)_{-n} (\log x)^t\\
&=& x^n \sum_{k=0}^t \r{n}! (t)_k s(-n,k) (\log x)^{t-k}.
\end{eqnarray*}
\end{enumerate}

Now we can give the association between the derivative and the
harmonic logarithms; the harmonic logarithms behave under derivation
exactly like the powers of $x$ once the ordinary factorial $n!$ is
replaced by the Roman factorial $\r{n}!$.
\begin{thm}\label{imported}
For all integers $n$, and vectors of integers with finite support
$\alpha$, we have
\begin{equation}\label{d}
\D \lambda_{n+1}^\alpha(x) = \r{n+1} \lambda_n^\alpha(x)
\end{equation}
where
$$\r{n}=\frac{\r{n}!}{\r{n-1}!} = \left\{ \begin{array}{ll}$$
                                n & \mbox{for $n\neq 0$, and}\\
                                1 & \mbox{for $n=0$.}
                                  \end{array}\right.$$
More generally, for any nonnegative integer $k$,
$$ \D ^k \lambda_{n+k}^\alpha = \frac{\r{n+k}!}{\r{n}!}
\lambda_n^\alpha (x).$$
\end{thm}

\proof By induction, it will suffice to demonstrate eq.~\ref{d}. The
following proof is a direct application of the recursion for the
Stirling numbers. (See \cite{stirl}.)
\begin{eqnarray*}
\lefteqn{\D\lambda_{n+1}^\alpha (x)}\\
 &=& \sum_{\rho}
\frac{s(-n-1,\rho_1)}{\r{-n-1}!} \left[ \prod_{i=1}^\infty e_{\rho_i -
\rho_{i+1}}
(\alpha_i, \alpha_i -1, \ldots, \alpha_i - \rho_i +1)\right]
\D \ell ^{(n+1),(\alpha-\rho)}\\
&=& \sum_{\rho}
\frac{s(-n-1,\rho_1)}{\r{-n-1}!} \left[ \prod_{i=1}^\infty e_{\rho_i -
\rho_{i+1}}
(\alpha_i, \alpha_i -1, \ldots, \alpha_i - \rho_i +1)\right]\\
&& \left[ (a+1) \ell ^{(n),(\alpha-\rho)}+ \sum_{k=1}^\infty (\alpha_k
- \rho_k ) \ell^{(n),(\alpha-\rho-(\overbrace{\scriptstyle
1,1,\ldots,1}^n)}) \right]\\
&=& \sum_{\rho}
\frac{(a+1)s(-n-1,\rho_1)+s(-n-1,\rho_1)}{\r{-n-1}!}\\
&& \left[ \prod_{i=1}^\infty e_{\rho_i - \rho_{i+1}}
(\alpha_i, \alpha_i -1, \ldots, \alpha_i - \rho_i +1)\right]
\ell^{(n),(\alpha-\rho)} \\
&=& \sum_{\rho}
\frac{s(-n,\rho_1)}{\r{-n-1}!}\left[ \prod_{i=1}^\infty e_{\rho_i - \rho_{i+1}}
(\alpha_i, \alpha_i -1, \ldots, \alpha_i - \rho_i +1)\right]
\ell^{(n),(\alpha-\rho)} \\
&=& \r{n+1} \lambda_n^\alpha (x).\Box
\end{eqnarray*}

Hence, if we denote by ${\cal L}^+$ the closure of the span of
the harmonic logarithm $\lambda_n^\alpha(x)$ with $\alpha\neq 0$, then
we notice that $\D$ restricted to ${\cal L}^+$ is a bijection.
In particular, the antiderivative is given by
$$ D^{-1} \lambda_n^\alpha (x) = \lambda_{n+1}^\alpha(x)/\r{n+1}.$$

We can thus characterize the harmonic logarithms by the
relationship.
$$ \lambda_n^\alpha(x) =
\r{n}!\D^{-n}\ell^{(0),\alpha}=\r{n}!\D^{-n}\lambda_0^{\alpha}.$$

Now, let us dissect the logarithmic algebra even more finely. Let
${\cal L}^\alpha$ denote the logarithmic series involving
$\lambda_n^\alpha(x)$ for all integers $n$ and a fixed $\alpha$.
Evidentally, ${\cal L}^+$ is the direct sum of the ${\cal L}^\alpha$
for $\alpha\neq(0)$.

\begin{thm} [Roman Modules]
All of the vector spaces  ${\cal L}^\alpha$ for $\alpha\neq(0)$ are
naturally isomorphic as differential vector spaces.
\end{thm}

\proof The isomorphism is given by the {\em skip operator}
$$\mbox{skip}_{\alpha\beta} \lambda_n^\beta (x)= \lambda_n^\alpha
(x).\Box$$

In view of the preceeding theory, most of our calculations can be made
in  ${\cal L}^{(1)}$ which is called the
{\em Roman module} after Steve Roman. The Roman module consists of
series of the form
$$ f(x)= \sum_{j=-\infty}^{-1} c_jx^j + \sum_{j=0}^n c_jx^j\left( \log
x - 1 - \frac{1}2 - \cdots \frac{1}j \right) = \sum_{j=-\infty}^n c_j
\lambda_j^{(1)}(x);$$
that is, there are a finite number of terms in which $\log x$ appears,
and an infinite sereies in inverse powers of $x$.
Thus, the Roman module is fairly convenient for our calculations. All
of our concrete examples will be drawn from the Roman module.

However, ${\cal L}^{(0)}$ is not isomorphic to the Roman module; it is
isomorphic to the algebra of polynomials ${\bf C}[x]$.

\subsection{Shift-Invariant Operators}

A linear operator $Q$ will be called {\em
shift-invariant} if not only does it commute with the shift operator
$E^a=\sum_{n=0}^\infty a^n\D^n/n! $ for all $a$, but also that it
commutes with the {\em skip operators} $\mbox{skip}_{\alpha\beta}$ from the
logarithmic algebra ${\cal L}$ to ${\cal L}^\alpha$
$$\mbox{skip}_{\alpha\beta}\lambda_n^\gamma(x) = \delta_{\beta\gamma}
\lambda_n^\alpha(x) $$
for $\beta\neq (0)$.

Any shift-invariant operator on ${\cal L}$ is obviously shift-invariant
(in the sense of eq.~\ref{eat}) when restricted to ${\cal L}^{(0)} = {\bf
C}[x]$. Conversely, any shift-invariant operator $Q$ on ${\bf C}[x]$ can
be expanded as a formal series in the derivative $Q(\D) $, and thus
$Q$ can be extended to all ${\cal L}$. However, are these the only
shift-invariant operators which act on ${\cal L}$? On ${\cal L}^+$?

\begin{thm} {\bf (Characterization of Shift-Invariant Operators)}

(1) The algebra of shift-invariant operators on ${\cal L}^+$ is
naturally isomorphic to the algebra of Laurent series\footnote{A {\em
Laurent series} is a formal series of the form
$$f(x)=\sum_{i=n}^\infty c_i x^i$$
where $n$ may be any constant positive or negative.}
in the derivative {\bf C}(\D).

(2) The algebra of shift-invariant operators on ${\cal L}$ is
naturally isomorphic to the algebra of formal power series
in the derivative {\bf C}[[\D]].
\end{thm}

\proof It will suffice to prove part one, since the only Laurent
series in the derivative which are well defined on ${\cal L}^{(0)}$ are
those which involve no negative powers of the derivative; that is,
members of the algebra of formal power series in the derivative.

Since the derivative commutes with $\mbox{skip}_{\alpha\beta}$ and
$E^a$, all Laurent series in the derivative are shift-invariant
operators on ${\cal L}^+$.

 Conversely, suppose that $Q$ is a shift-invariant operator. $Q$ is
determined by its actions on the harmonic logarithms of order $\alpha$
for any particular $\alpha\neq (0)$.
$$Q \lambda_n^\alpha (x) = \sum_{k=-\infty}^n c_{nk} \lambda_k^\alpha
(x).$$
Next, since $Q$ commutes with $E^a$ for all $a$, it also commutes with
$D^n$ for all $n$. Thus,
\begin{eqnarray*}
 \sum_{j=-\infty}^k c_{kj}
\frac{\r{k}!}{\r{k-j}!}\lambda_{k-j}^\alpha(x)
&=& \D^nQ\lambda_k^\alpha(x)\\
& =& Q\D^n \lambda_k^\alpha(x) \\
&=& \frac{\r{k}!}{\r{k-n}!}\sum_{j=-\infty}^{k-n} c_{n-k,j}
\lambda_j^\alpha (x).
\end{eqnarray*}
Equating coefficients of $\lambda_n^\alpha(x)$, and setting $a=b$, we have
$$ c_{0j} = \rc{j}{k}  c_{k,j+k} $$
where the {\em Roman coefficient} $\rc{j}{k}$ is defined to be
$\r{j}!/\r{k}!\r{j-k}!$.

Hence, $Q$ is determined by the $c_{0j}$, and
therefore equals the Laurent series
$\sum_{j=-\infty}^{m}c_{0j}\D^j/\r{j}!$ where $m$ is the largest $j$
such that $c_{0j}\neq 0. \Box$

Since the algebra of Laurent series is a field, we immediately derive
that all shift-invariant operators on ${\cal L}^+$ are invertible.
Thus,

\begin{cor}[Differential Equations]\label{sub}

(1) All linear differential equations with constant coefficients on
${\cal L}^+$ have a {\em 
unique} solution. That is, any equation of the form
$$f(\D) p(x)= q(x)$$
(where $q(x)\in {\cal L}^+$)
has a unique solution $p(x)\in {\cal L}^+$.

(2) All linear differential equations on ${\cal L}$ have solutions.
Of these, one particular solution can be
naturally chosen as the canonical one.
\end{cor}

\proof{We have just proven 1. The proof for 2  (see \cite[\S
2.4.4A]{ILA1}) relies on the
 use of the projection maps $\mbox{skip}_\alpha\beta.\Box$

\subsection{Binomial Theorem}

Any shift-invariant operator acts on each level ${\cal L}^\alpha$ in
the same way. For example, consider the logarithmic version of the
binomial theorem.
\begin{equation}\label{log-bin}
 E^a \lambda_n^\alpha (x) = \sum_{k=0}^\infty a^k \D^k
\lambda_n^\alpha(x)/k! = \sum_{k=0}^\infty \rc{n}{k} a^k
\lambda_{n-k}^\alpha (x).
\end{equation}
Several special cases are of particular interest. For $\alpha=(0)$, we
have the standard binomial theorem. For $\alpha=(1)$, $n$ positive,
and $x=1$, we have \cite[p. 97]{thesis}
\begin{eqnarray*}
\lefteqn{(1+a)^n\log(1+a)}\\
 &=& ((1+a)^n -1)\left(1+\frac{1}2 +\cdots +
\frac{1}n\right) -na\left(1+\frac{1}2+\cdots+\frac{1}{n-1}\right)\\
&& -{n \choose 2} a^2 \left(1+\frac{1}2+\cdots+\frac{1}{n-2}\right) -\cdots
-\frac{3}{2}{n \choose n-2} a^{n-2} -n a^{n-1} \\
&& +a^{n+1}\rc{n}{-1} + a^{n+2}\rc{n}{-2}+\cdots.
\end{eqnarray*}
Similarly, for $\alpha=(\overbrace{0,0,\ldots,0}^{k-1},1)$ and $n=0$,
\begin{eqnarray*}
\lefteqn{\overbrace{\log \log \cdots \log}^{n}(x+a)}\\
&=& \overbrace{\log \log \cdots \log}^{n}(x)
 +\sum_{\rho_1 \geq \rho_2 \geq \cdots \geq \rho_k = 1}
(-1)^{\rho_1+\cdots+\rho_{k-2}+2\rho_{k-1}+1}\\
&& \times \left[ \prod_{j=1}^{k-2} e_{\rho_j - \rho_{j+1}}(1,\ldots,
\rho_j -1)\right]
\sum_{n=0}^\infty
z^n(n-1)!s(-n,\rho_1)\ell^{(-n,\rho_1,\ldots,\rho_{k-1})}
\end{eqnarray*}

This framework immediately gives
rise to logarithmic versions of Taylor's theorem and the expansion
theorem in terms of the logarithmic generalization of evaluation at
zero or augmentation. The {\em augmentation}\/\footnote{The restriction
of the augmentation of order $(0)$ to ${\cal L}^{(0)}$ is the
evaluation at zero mentioned in eq.~\ref{augment}. The other
augmentations can be profitably be rewritten in terms of the residue
at zero.} of order $\alpha$ is defined to be the linear functional
$\aug{}_{\alpha}$ from the logarithmic algebra ${\cal L}$ to the
complex number {\bf C} such that
$$ \aug{\lambda_n^{\beta}(x)}_\alpha = \delta_{\alpha\beta}
\delta_{n0}.$$

\begin{thm}[Logarithmic Taylor's Theorem] \label{tailor}

Let $p(x)\in {\cal L}^+$. Then we have the following expansion of
$p(x)$ in terms of the harmonic logarithms.
\begin{equation}\label{holds}
p(x)=\sum_{\alpha\neq (0)} \sum_{n=\infty}^\infty
\frac{\aug{\D^np(x)}_\alpha}{\r{n}!} \lambda_n^\alpha(x).
\end{equation}
\end{thm}

\proof It suffices to show that eq.~\ref{holds} holds for a basis of
${\cal L}^+$. However, Theorem~\ref{imported} states that eq.~\ref{holds}
holds for the basis of harmonic logarithms.$\Box$

\begin{thm}[Logarithmic Expansion Theorem] \label{let}

(1) Let $Q(\D)$ be a shift-invariant operator on ${\cal L}^+$. Then we
have the following expansion of $Q(\D)$ in terms of the powers of the
derivative
$$ Q(\D) = \sum_{n=-\infty}^\infty
\frac{\aug{Q(\D)\lambda_n^\alpha(x)}_\alpha}{\r{n}!} \D^n $$
where $\alpha\neq (0)$.

(2) Let $Q(\D)$ be a shift-invariant operator on ${\cal L}$. Then we
have the following expansion of $Q(\D)$ in terms of the powers of the
derivative
$$ Q(\D) = \sum_{n=0}^\infty
\frac{\aug{Q(\D)\lambda_n^\alpha(x)}_\alpha}{n!} \D^n $$
where $\alpha$ is any vector of integers with finite support.
\end{thm}

\proof Theorem is immediately for $Q(\D)=\D^n$ from Theorem
\ref{imported}.$\Box$

We now see that the sequence of harmonic logariths
$\lambda_n^\alpha(x)$ is in some sense (to
be made clear later) of ``binomial type,'' and that in this sense the
harmonic logarithms are associated with the derivative operator.

\subsection{The  Logarithmic Lower
Factorial}\label{eggs}

We can now extend to the logarithmic domain the analogy between the
derivative operator $\D$ and the forward difference operator $\Delta$
which has been noted for polynomials.

Since ${\cal L}^(0)$ is the algebra of polynomials, one would rightly
expect that the logarithmic version $p_n^\alpha(x)=(x)_n^\alpha$ of
the lower factorial function $p_n(x)=(x)_n$ would include the
subsequence
$$p_n^{(0)}(x)=(x)_n^{(0)}=\left\{
\begin{array}{ll}
(x)_n = x(x-1)\cdots (x-n+1) & \mbox{for $n$ nonnegative, and}\\
0 &\mbox{for $n$ negative.}
\end{array}\right.$$

It has been known for a long time that, setting for $n$ negative,
$$p_n(x)= (x)_n = 1/(x+1)(x+2)\cdots(x-n),$$
one has
$$ \Delta p_n(x)= n(x)_{n-1}.$$
Given this preliminary information,\footnote{This information is only
preliminary since without the residual theorem
(Theorem~\ref{residual}) to pinpoint the exact
value of $p_{-1}(x)$, the entire sequence might be off by a factor
equal to a shift-invariant operator of degree zero
$a_0+a_1\D+a_2D^2+\cdots.$} we  guess that when $n$ is
negative, $1/(x+1)\cdots(x-n)$ is to the forward difference operator
as $x^{n}$ is to the derivative. Thus, we write
$$p_n^{(1)}(x)=(x)_n^(1)=1/(x+1)(x+2)\cdots(x+n)$$
for $n$ negative.

This leave open the question of what is the forward difference version
of the logarithm $\lambda_0^{(1)}(x)=\log x$ and of the harmonic
logarithms $$\lambda_n^{(1)}(x)=x^n\left(\log x-
1-\frac{1}2-\cdots-\frac{1}n\right) $$  for $n$ positive.
In other words, what is the formal solution $p_0^{(1)}(x)=(x)_0^{(1)}$
of the difference equation $\Delta p(x)=1/(x+1)$? This question was
first formulated by Gauss; in
the present context, it can be dealt with very easily. By
Corollary~\ref{sub}, all such difference equations
have a unique solutions.

In fact, we can easily calculate $(x)_0^{(1)}$. Since, the {\em Bernoulli
numbers} $B_n$ are defined as the coefficients of $\Delta^{-1}$:
$$\frac{\D}{e^{\D} -1} = \sum_{k=0}^\infty B_k\D^k/k!,$$
we have
$$p_0^{(1)}(x)=(x)_0^{(1)}=\log(x+1)+B_1/(x+1)-B_2/2(x+1)^2
+B_3/3(x+1)^3 -\cdots. $$
Note that $(x-1)_0^{(1)}=\psi(x)$ is the {\em Gauss psi function}. In
this domain, we verify that
$$\psi(x)=\Gamma'(x)/\Gamma(x)$$
as follows. Begin with the fundamental identity for the Gamma
function,
$$\Gamma (x+1) = x\Gamma(x).$$
Now, take the derivative of both sides.
$$ \Gamma'(x+1) = \Gamma (x) + x\Gamma (x).$$
Next, divide both sides by $\Gamma(x+1)=x\Gamma(x).$
$$ \frac{\Gamma'(x+1)}{\Gamma(x+1)} = \frac{1}{x} +
\frac{\Gamma'(x)}{\Gamma(x)}.$$
Finally, note that
$$ \frac{1}{x} = \Delta \frac{\Gamma'(x)}{\Gamma(x)}.$$
Thus, $\Gamma'(x)/\Gamma(x)$ equals $(x)_0^{(1)}$ possibly up to a
constant.

For the logarithmic lower factorial of order (1) and positive degree $n$, we
simply employ the higher order Bernoulli numbers
$$p_n^{(1)}(x)=(x)_n^{(1)}=\sum_{k=1}^\infty B_{k,n+1} \rc{n}{k}
\lambda_{n-k}^\alpha(x+1)$$
where
$$\left(\frac{\D}{e^{\D} -1}\right)^n = \sum_{k=0}^\infty B_{kn}\D^k/k!.$$

Finally, to compute the logarithmic lower factorials of order
$\alpha\neq (1)$, merely use the projection maps.
$$p_n^\alpha(x)= (x)_n^\alpha = \mbox{skip}_{\alpha,(1)}
(x)_n^{(1)}.$$

To carry the analogy to the end, we note that by eq.~\ref{taylor} and
part 2 of Theorem~\ref{let}, we still have
$$ E^a = \sum_{n=0}^\infty (a)_n \Delta^n/n!.$$
From this, we have the logarithmic version of VanderMonde's identity
$$ (x+a)_n^{\alpha} = \sum_{k=0}^\infty \rc{n}{k} (a)_k
(x)_{n-k}^\alpha$$
for all integers $n$. Similarly, we have the logarithmic analog of the
Newton's formula
$$ f(x) = \sum_{\alpha\neq (0)}\sum_{n=-\infty}^\infty \aug{\Delta^n
f(x)}_\alpha (x)_n^\alpha/n!.$$
We note that this formula constitutes a substantial extension of
Noerton's formula. For example, if $f(x)= 1/x$, then we have for $n$
negative\footnote{Obviously, for $n$ nonnegative, $\aug{\Delta^n
\left(\frac{1}{x}\right) }_{(1)}=0$.}
\begin{eqnarray*}
\Delta^n \frac{1}x &=& \aug{E^{-1} \Delta^n (x)_{-1}^{(1)}}_{(1)} \\
&=& \frac{\aug{E^{-1} (x)_{-n-1}^{(1)}}_{(1)}}{(-n-1)!}   \\
&=& \frac{\aug{E^{-1} (x)_{-n-1}^{(0)}}_{(0)}}{(-n-1)!}   \\
&=& \left[ (x)_{-n-1}^{(0)}/(-n-1)! \right]_{x=-1}  \\
&=& (-1)(-2)\cdots(n)/(-n-1)!\\
&=& (-1)^{n+1}n .
\end{eqnarray*}
Thus,
\begin{eqnarray}
\frac{1}{x} &=& \sum_{n=-\infty}^{-1} \frac{(-1)^{n+1}n}{\r{n}!}
(x)_n^{(1)}\nonumber\\
&=& \sum_{n=0}^\infty \frac{n!}{(x+1)(x+2)\cdots(x+n)}.\label{golden}
\end{eqnarray}

\subsection{Logarithmic Sequences of Binomial Type}
As we did with polynomials, we too can now extend the interplay
between the derivative $\D$ and the forward difference operator into a
more general framework involving all delta operators.\footnote{As
before, a delta operator is a delta series in the derivative; that is
a shift-invariant operator whose kernel is the set of constants---in
this case ${\bf R}$.} Instead of sequences of polynomials $p_n(x)$, we
have {\em logarithmic sequences} $p_n^\alpha (x)$ indexed by an integer $n$
and a vector with finite support of integers $\alpha$. Again we have a
requirement that $p_n^\alpha (x)$ be of {\em degree} $n$, and also a
new requirement that it be of {\em order} $\alpha$; that is to say,
$$p_n^\alpha (x) = \sum_{-\infty}^n c_k \lambda_k^\alpha (x) $$
with $c_n \neq 0$.
In particular, $p_n^{(0)}(x)$ is a sequence of polynomials: one of each
degree with $p_n^{(0)}(x)=0$ for $n$ negative. More generally, the
$p_n^\alpha (x)$ each belong to ${\cal L}^\alpha$ and in fact form a
basis for it.

The only two examples of logarithmic sequences which we have seen so far
are the harmonic logarithm $\lambda_n^\alpha(x)$ and the logarithmic
lower factorial $(x)_n^\alpha$. However, there are clearly many more;
every sequence of polynomials can be extended in many ways into a
logarithmic sequence. However, we are particularly
interested in sequences of polynomials of binomial type. We seek a
natural definition of a logarithmic sequence of
binomial type such that every sequence of polynomials has one unique
such extension. The following theorem is exactly what we are seeking.

\begin{thm} {\bf (Logarithmic Sequences  of Binomial Type)}
\label{nome}
Let $p_n(x)$ be a sequence of polynomials of binomial
type. Then there exists a unique logarithmic sequence $p_n^\alpha (x)$
called the {\em logarithmic sequence of binomial type}
which
\begin{enumerate}
\index extends $p_n(x)$
\begin{equation}\label{a}
p_n^{(0)}(x) = p_n(x) \mbox{ for all integers $n$,}
\end{equation}
\index obeys an analog of the binomial theorem
\begin{equation}\label{b}
p_n^\alpha(x+a) = \sum_{k=0}^\infty \rc{n}{k} p_k(a) p_{n-k}^\alpha(x)
\end{equation}
for all integers $n$ and constants $a$,
\index and is  invariant under the skip operator
$\mbox{skip}_{\alpha\beta}$ for $\beta\neq(0)$; that is to say,
\begin{equation}\label{skipped}
 \mbox{skip}_{\alpha\beta}p_n^\gamma(x) = \delta_{\beta\gamma}
p_n^\alpha (x).
\end{equation}
\end{enumerate}
\end{thm}

\proof {\bf Existence} Let $f(\D)$ be the delta operator associated with
$p_n(x)$ and let $g(\D)$ be its compositional inverse.\footnote{They
are compositional inverses in the sense of composition of formal power
series in the derivative. That is, $f(g(\D))=\D=g(f(\D))$} Now define
$p_n^\alpha(x)$ to be what is known commonly as the {\em conjugate
logarithmic sequence} for $g(\D)$.
$$p_n^\alpha(x)=\sum_{k=-\infty}^n \frac{\aug{g(\D)^k
\lambda_n^\alpha(x)}_\alpha}{\r{n}!} \lambda_k^\alpha(x).$$
After briefly noting that $p_n^\alpha(x)$ is a well defined logarithmic
sequence invariant under the skip operators, it will now suffice to
prove it obeys eqs.~\ref{a} and
\ref{b}. This will follow immediately from the following two lemmas.

\begin{lem}[Associated Logarithmic Sequences] Let $p_n^\alpha(x)$ be the
conjugate logarithmic sequence for
the delta operator $g(\D)$ whose compositional inverse is $f(\D)$. Then
$p_n^\alpha(x)$ is the {\em associated logarithmic sequence} for $f(\D)$ in
the sense that
\begin{eqnarray}
\aug{p_n^\alpha(x)}_\alpha &=& \delta_{n0} \label{one}\\
f(\D) p_n^\alpha (x) &=& \r{n} p_{n-1}^{\alpha}(x). \label{two}
\end{eqnarray}
\end{lem}

\proof Eq.~\ref{one} follows directly from the definition of a
conjugate logarithmic sequence and the augmentation,
$$\aug{p_n^\alpha (x)}_\alpha = \aug{\lambda_n^\alpha (x)}/\r{0}! =
\delta_{n0} . $$

For eq.~\ref{two}, we must first write $f(\D)$ in terms of $\D^n$:
$f(\D) = \sum_{k=0}^\infty a_k \D ^k$. Then we proceed as follow
\begin{eqnarray*}
f(\D) p_n^\alpha (x) &=& \sum_{k,j\geq 0} a_j \frac{\aug{g(\D)^k
\lambda_n^\alpha(x)}_\alpha}{\r{k-j}!} \lambda_{k-j}^\alpha(x)\\
 &=& \sum_{k,j\geq 0}  \frac{\aug{a_j g(\D)^{k+j}
\lambda_n^\alpha(x)}_\alpha}{\r{k}!} \lambda_{k}^\alpha(x)\\
 &=& \sum_{k= 0}^\infty  \frac{\aug{\left(\sum_{j=0}^\infty a_j
g(\D)^{j} \right) g(\D)^k
\lambda_n^\alpha(x)}_\alpha}{\r{k}!} \lambda_{k}^\alpha(x)\\
 &=& \sum_{k= 0}^\infty  \frac{\aug{ \D g(\D)^k
\lambda_n^\alpha(x)}_\alpha}{\r{k}!} \lambda_{k}^\alpha(x)\\
 &=& \r{n} \sum_{k= 0}^\infty  \frac{\aug{ g(\D)^k
\lambda_{n-1}^\alpha(x)}_\alpha}{\r{k}!} \lambda_{k}^\alpha(x)\\
&=& \r{n} p_{n-1}^\alpha(x).\Box
\end{eqnarray*}

Note that in the case $\alpha=(0)$, $p_n^{(0)}(0)=\delta_{n0}$ and
$f(\D) p_n^{(0)}(x)=np_{n-1}^{(0)}(x)$, so that $p_n^{(0)}(x)$ is
associated to the delta operator $f(\D)$. Hence, we have eq.~\ref{a}.

\begin{lem}\label{lemma}
Let $p_n^\alpha (x)$ be the associated logarithmic sequence for the delta
operator $f(\D)$. Then for all shift-invariant operators $h(\D)$ we
have
$$ h(\D) p_n^\alpha (x) = \sum_{k=-\infty}^\infty \rc{n}{k}
\aug{h(\D)p_k^\alpha(x)}_\alpha p_{n-k}^\alpha(x).$$
\end{lem}

\proof It suffices to prove the lemma for $h(\D)=f(\D)^j$ since the
powers of the delta operator $f(\D)$ form a basis for the space of
shift-invariant operators. However, by the definition of an associated
logarithmic sequence, we immediately $\aug{f(\D)^j
p_k^\alpha(x)}=\delta_{jk}\r{j}!$. Thus,
$$ f(\D)^j p_n^\alpha (x) =\frac{\r{n}!}{\r{n-j}!}
\sum_{k=-\infty}^\infty \rc{n}{k}
\aug{f(\D)^jp_k^\alpha(x)}_\alpha p_{n-k}^\alpha(x).\Box$$

Finally, eq.~\ref{b} follows from the lemma with $h(\D)=E^a$.

{\bf Uniqueness} Now it remains only to show the uniquness of
logarithmic sequences   of binomial type. We will do so
in two lemmas. The first shows that all logarithmic sequences of binomial
type are associated logarithmic sequences, and the other  shows that every
delta operator has a unique associated logarithmic sequence. These two
lemmas will complete the proof of Theorem~\ref{nome}.

\begin{lem}
Let $p_n^\alpha(x)$ be a logarithmic sequence of binomial type; that is,
let it obey eq.~\ref{b}. Then there is a unique delta operator $f(\D)$
such that $p_n^\alpha(x)$ is associated with $f(\D)$. That is,
$p_n^\alpha(x)$ obeys eqs.~\ref{one} and~\ref{two}.
\end{lem}

\proof Eq.~\ref{one} immediately follows from
Proposition~\ref{zero}. Now, we {\em define} a linear operator $Q$ by
$$Qp_n^\alpha(x) = \r{n} p_{n-1}^\alpha(x).$$
Eq.~\ref{two} is equivalent to the proposition that $Q$ is a delta
operator. Since $Q$ obviously commutes with
$\mbox{skip}_{\alpha\beta}$ for $\beta\neq(0)$, and the kernel of $Q$
is clearly the set of real constants {\bf R}, it remains only to show
that $Q$ commutes with the shift operator $E^a$.
\begin{eqnarray*}
QE^a p_n^\alpha(x) &=&  \sum_{k=0}^\infty \rc{n}{k} p_k(a)
Qp_{n-k}^\alpha(x)\\
&=&  \sum_{k=0}^\infty \rc{n}{k} \r{n-k} p_k(a) p_{n-k-1}^\alpha(x)\\
&=&  \sum_{k=0}^\infty \rc{n-1}{k} \r{n} p_k(a) p_{n-k-1}^\alpha(x)\\
&=&  E^a Qp_n^\alpha(x).\Box
\end{eqnarray*}

\begin{lem}\label{eunuch} Let $f(\D)$ be a delta operator. Then there
is exactly one
logarithmic sequence $p_n^{\alpha}(x)$ associated
with  $f(\D)$.
\end{lem}

\proof By Taylor's theorem (Theorem~\ref{tailor}), $p_n^\alpha(x)$ is
determined by its augmentations $a_k=\aug{\D^k p_n^\alpha(x)}_\alpha$.
However, since $f(\D)^k$ is a basis for the space of shift-invariant
operators, the augmentations $a_k$ are determined by the augmentations
$b_k=\aug{f(\D)^k p_n^\alpha(x)}_\alpha$ which we know are equal to
$\r{n}!\delta_{nk}.\Box$ QED

We also derive---in the same manner as Lemma
\ref{lemma}---generalizations of the logarithmic Taylor's theorem and
expansion theorem.
\begin{thm}[Logarithmic Taylor's Theorem] \label{tail}
Let $p_n^\alpha(x)$ be the logarithmic sequence associated with the delta
operator $f(\D).$ Then any logarithmic series $p(x)\in {\cal L}^+$ can
be  expanded
$$ p(x) = \sum_{\alpha\neq (0)} \sum_{n=-\infty}^{N_\alpha}
\frac{\aug{f(\D)^np(x)}_\alpha}{\r{n}!} p_n^\alpha (x).$$
\end{thm}
\begin{thm}[Logarithmic Expansion Theorem]
Let $p_n^\alpha(x)$ be the logarithmic sequence associated with the delta
operator $f(\D),$ and let $\alpha\neq(0)$. Then any shift-invariant
operator $g(\D)$ can be  expanded
$$ g(x) =  \sum_{n=-\infty}^\infty
\frac{\aug{g(\D)p_n^\alpha(x)}_\alpha}{\r{n}!} f(\D)^n.$$
\end{thm}

\subsection{Explicit Formulas}

We note that if any term of the logarithmic sequence  of binomial type
$p_n^{(1)}(x)$ associated with the
delta operator $f(\D)$ is known a priori,
then all of the other terms can be computed in terms of it by
appropriate use of the various powers of $f(\D)$ and the projections
$\mbox{skip}_{\alpha,(1)}$. In particular, we only need to determine
the term $p_{-1}^{(1)}(x)$ which is called the {\em residual term},
and then
$$p_n^\alpha (x) =
\mbox{skip}_{\alpha,(1)}f(\D)^{-n-1}p_{-1}^{(1)}(x)/\r{n}!.$$

Thus, the following result is extemely fundamental,
\begin{thm}[Residual Term]\label{residual}
Let $p_n^\alpha(x)$ be the logarithmic sequence associated with the delta
operator $f(\D),$ and let $\alpha\neq(0)$. Then\footnote{The operator
$f'(\D)$ called the {\em Pincherle derivative} of $f(\D)$ represents
the operator $g(\D)$ where
$g(t)=\frac{df(t)}{dt}$. Thus, if $f(\D)=\D^n$, then
$f'(\D)=n\D^{n-1}.$ Another characterization of the Pincherle
derivative will be given later. (Theorem~\ref{pinch})}
$$ p_{-1}^{(1)}(x) = f'(\D) \frac{1}{x}. $$
\end{thm}

\proof We will actually prove a seemingly stronger result
\begin{equation}\label{seams}
p_n^\alpha(x) = \r{n}! f'(\D) f(\D)^{-1-n} \lambda_{-1}^\alpha(x).
\end{equation}
Let $q_n^\alpha(x)$ denote the right hand side of eq.~\ref{seams}.
Clearly, $q_n^\alpha(x)$ is a logarithmic sequence, so by Lemma
\ref{eunuch}, it will suffice to
show that $q_n^\alpha(x)$ is associated with the delta operator
$f(\D)$. Of the two defining properties, eq.~\ref{two} is trivial to
demonstrate.
$$ f(\D)q_n^\alpha(x) = \r{n} \r{n-1}! f'(\D) f(\D)^{-n}
\lambda_{-1}^\alpha(x) = \r{n} q_{n-1}^\alpha(x).$$
To prove eq.~\ref{one}, we first consider the case $n\neq 0$.
\begin{eqnarray*}
\aug{p_n^\alpha(x)}_\alpha &=& \aug{\r{n}! f'(\D) f(\D)^{-1-n}
\lambda_{-1}^\alpha(x)}_\alpha \\
&=& \frac{\r{n}!}{-n} \aug{ \left(f(\D)^{-n}\right)'
\lambda_{-1}^\alpha(x)}_\alpha
\end{eqnarray*}
However, $\aug{ \left(f(\D)^{-n}\right)'
\lambda_{-1}^\alpha(x)}_\alpha $ is given by the coefficient of
$\D^{-1}$ in $\left(f(\D)^{-n}\right)'$; however, this coefficient is
always zero.\footnote{This theorem will be seen to be equivalent to
the Lagrange inversion formula (Theorem~\ref{Large}). The observation here
that the residue
of the derivative of a formal power series is zero is the key step in
all known proofs of the Lagrange inversion formula.} Thus,
$\aug{p_n^\alpha(x)}_\alpha = 0$ for $n\neq 0$.

Whereas for $n=0$, we have
$$\aug{p_0^\alpha(x)}_\alpha = \aug{ f'(\D) f(\D)^{-1}
\lambda_{-1}^\alpha(x)}_\alpha .$$
However, $f'(\D)f(\D)^-1$ is equal to $\D^{-1}$ plus terms of higher
degree, so
$$\aug{p_0^\alpha(x)}_\alpha = 1.\Box$$

Theorem~\ref{residual} allows us to derive a series of formulas---as
explicit as possible---for the computation of sequences of polynomials
of binomial type or logarithmic sequence  of binomial
type associated with a delta operator $f(\D)$.

However, we require one last definition. We define on the logarithmic
algebra the {\em Roman shift}\/\footnote{Named after Steve Roman.
Sometimes also called the {\em Standard Roman Shift} or the {\em Roman
Shift for $\D$} or the {\em Roman Shift for $\lambda_n^\alpha(x)$.}}
 as the linear operator $\sigma$ such that for all $\alpha$
$$\sigma \lambda_n^\alpha (x) = \left\{
\begin{array}{ll}
\lambda_{n+1}^\alpha (x) & \mbox{for $n\neq -1$, and}\\
0& \mbox{for $n=-1$.}
\end{array}\right.$$
The Roman shift is the logarithmic generalization of multiplication by
$x$. In particular, note that for $n\neq -1$,
$$\sigma x^n = x^{n+1}. $$

We observe the following remarkable facts about the Roman shift.
\begin{prop}[Pincherle Derivative] \label{pinch}
For all logarithmic series $p(x)\in \cal L$,
\begin{equation}\label{vain}
(\D \sigma - \sigma \D ) p(x) = p(x);
\end{equation}
in other words,
$$\D \sigma - \sigma \D = I.$$
More generally, for all shift-invariant operators $f(\D)$, we have the
operator identity
\begin{equation}\label{sigmund}
 f(\D) \sigma - \sigma f(\D) = f'(\D).
\end{equation}
\end{prop}

\proof It will suffice to consider the case $f(\D)=\D^k$ and
$\lambda_n^\alpha(x)$. When $k=0$, eq.~\ref{sigmund} is a triviality.
Now, consider the case in which $n$ is neither $k-1$ nor $-1$:
\begin{eqnarray*}
(\D^k \sigma -\sigma \D^k) \lambda_n^\alpha (x) &=&
\D ^k \lambda_{n+1}^\alpha (x) -\frac{\r{n}!}{\r{n-k}!}\sigma
\lambda_{n-k}^\alpha(x) \\
&=& \left(\frac{\r{n+1}!}{\r{n-k+1}!} - \frac{\r{n}!}{\r{n-k}!}\right)
\lambda_{n-k+1}^\alpha(x) \\
&=& (\r{n+1}-\r{n-k+1}) \frac{\r{n}!}{\r{n-k+1}!}
\lambda_{n-k+1}^\alpha(x) \\
&=& k \D^{k-1}\lambda_{n}^\alpha(x).
\end{eqnarray*}
Next, consider the case in which $n= -1$ but $k\neq 0$:
\begin{eqnarray*}
(\D^k \sigma -\sigma \D^k) \lambda_{-1}^\alpha (x) &=& -\sigma \D^k
\lambda_{-1}^\alpha (x)\\
&=& -\frac{1}{\r{-k-1}!}\sigma
\lambda_{-k-1}^\alpha(x) \\
&=& -\r{-k}\frac{1}{\r{-k}!}
\lambda_{-k}^\alpha(x) \\
&=& k \D^{k-1}\lambda_{-1}^\alpha(x) .
\end{eqnarray*}
Finally, consider the case in which $n=k-1$:
\begin{eqnarray*}
(\D^k \sigma -\sigma \D^k) \lambda_n^\alpha (x) &=&
\D ^k \lambda_{n+1}^\alpha (x) -\frac{\r{n}!}{\r{n-k}!}\sigma
\lambda_{-1}^\alpha(x) \\
&=& \r{n+1}! \lambda_{0}^\alpha(x) \\
&=& \r{n+1} \D^{n} \lambda_{n}^\alpha(x) \\
&=& k \D^{k-1}\lambda_{n}^\alpha(x).\Box
\end{eqnarray*}

Thus, we see that the Roman shift gives a non-trivial extension of the
classical commutation relation $\D x-x\D= I$ of quantum mechanics in
the logarithmic algebra. In effect, eq.~\ref{vain} is really quite
incredible, for Van Neuman has shown \cite{VN} that for
polynomials\footnote{Since then others have shown that this result
holds for much more general sorts of functions although obviously not
for logarithmic series.} the only operators $(\theta, \phi)$ obeying
the commutation relation.
$$\theta \phi - \phi \theta = I$$
are $\theta= \D$ and $\phi = x$. However, by eq.~\ref{vain}, we see
that for logarithmic series, $\theta = \D$ and $\phi = \sigma$ is yet
another solution.


\begin{thm}[Rorigues]
Let $p(x)$ be a finite linear combination of the harmonic logarithms.
Then for all real numbers $a$, we have the following formal identity
$$ (\exp(-a\sigma)\D \exp(a\sigma)) p(x) = (\D-aI)p(x).\Box$$
\end{thm}

We can now state our main result.
\begin{thm}\label{mane}
Let $f(\D)=\D g(\D)$ be the delta operator associated both with  the
logarithmic 
sequence $p_n^\alpha (x)$, and with the sequence
of polynomials $p_n(x)=p_n^{(0)}(x)$. Then we have,
\begin{enumerate}
\item {\bf (Transfer Formula)} For all $n$,
$$ p_n^\alpha (x) = f'(\D) g(\D)^{-n-1} \lambda_n^\alpha (x),$$
\item {\bf (Shift Formula)} For $n$ neither 0 nor $-1$,
$$ p_n^\alpha (x) = \sigma g(\D)^{-n} \lambda_{n-1}^\alpha (x),$$
and
\item {\bf (Recurrence Formula)} For $n\neq -1$,
$$ p_n^\alpha (x) = \sigma f'(\D)^{-1} p_{n-1}^\alpha (x).$$
\end{enumerate}

In particular, for the polynomial sequence $p_n(x)$, we have
\begin{enumerate}
\item {\bf (Transfer Formula)}
$$ p_n(x) = f'(\D) g(\D)^{-n-1} x^{n},$$
\item {\bf (Shift Formula)}
$$ p_n(x) = x g(\D)^{-n} x^{n-1},$$
and
\item {\bf (Recurrence Formula)}
$$ p_n(x) = x f'(\D)^{-1} p_{n-1}(x).$$
\end{enumerate}
\end{thm}

\proof (1) From eq.~\ref{seams}, we have for $\alpha\neq(0)$
$$ p_n^\alpha(x) = \r{n}! f'(\D) g(\D)^{n+1} \D^{-1-n}
\lambda_{-1}^\alpha(x) =  f'(\D) g(\D)^{n+1} \lambda_n^\alpha (x).$$
For $\alpha=(0)$, we derive the result from the case $\alpha=(1)$ as
follows\footnote{This is then the best possible example of a purely
logarithmic proof of a result involving polynomials. It is derived
from a formula (eq.~\ref{seams}) which does not even make sense for
polynomials.}
\begin{eqnarray*}
p_n(x) &=& \mbox{skip}_{(0),(1)} p_n^{(1)}(x) \\ 
 &=& \mbox{skip}_{(0),(1)} f'(\D) g(\D) ^{n+1} \lambda_n^{(1)}(x) \\
 &=& f'(\D) g(\D) ^{n+1} \mbox{skip}_{(0),(1)} \lambda_n^{(1)}(x) \\
 &=& f'(\D) g(\D) ^{n+1} x^n.
\end{eqnarray*}

(2) For this part, we need the following routine calculation
$$f'(\D) g(\D) ^{n+1} = g(\D)^{n} - \frac{1}{n}(g(\D)^{n})'\D .$$
Thus, by part (1), we have
$$p_n^\alpha (x) = g(\D)^{n}\lambda_n^\alpha(x) - (g(\D)^{n})'
\lambda_{n-1}^\alpha(x).$$
Finally, by Proposition~\ref{pinch}, we have
\begin{eqnarray*}
p_n^\alpha (x) &=& g(\D)^{n}\lambda_n^\alpha(x) - g(\D)^{n}\sigma
\lambda_{n-1}^\alpha(x)+ \sigma g(\D)^{n}
\lambda_{n-1}^\alpha(x) \\
&=& g(\D)^{n}\lambda_n^\alpha(x) - g(\D)^{n}
\lambda_{n}^\alpha(x)+ \sigma g(\D)^{n}
\lambda_{n-1}^\alpha(x) \\
&=& \sigma g(\D)^{n} \lambda_{n-1}^\alpha(x).
\end{eqnarray*}

(3) This follows immediately from the substitute of the Transfer
Formula
$$ \lambda_{n-1}^\alpha (x) =  f'(\D)^{-1} g(\D)^{n}p_{n-1}^\alpha
(x)$$
in the Shift Formula $$ p_n^\alpha (x) = \sigma g(\D)^{-n}
\lambda_{n-1}^\alpha (x).\Box$$

Actually, Theorem~\ref{mane} can be generalized to cover the case in
which the role of the harmonic logarithm $\lambda_n^\alpha (x)$ can be
played by any logarithmic sequence of binomial type $p_n^\alpha (x)$.
\begin{cor}\label{cute}
Let $f(\D)$ and $g(\D)$ be delta operators associated with the
the logarithmic sequences $p_n^\alpha(x)$ and $q_n^\alpha (x)$, and
the sequences of polynomials $p_n(x)=p_n^{(0)}(x)$ and
$q_n(x)=q_n^{(0)}(x)$.
Then
\begin{enumerate}
\item We have
$$ p_n^\alpha (x) = \frac{f'(\D)g(\D)^{n+1}}{g'(\D)f(\D)^{n+1}}
q_n^\alpha (x) $$
for all $\alpha$ and $n$, and in particular for $\alpha = (0) $ we
have
$$ p_n(x) = \frac{f'(\D)g(\D)^{n+1}}{g'(\D)f(\D)^{n+1}} q_n (x) ,$$
\item And for the polynomial sequences only,\footnote{The Roman shift
is not in general invertible.} we have
$$ p_n(x) = x\left( g(\D)^n f(\D)^{-n} \left[ x^{-1} q_n (x) \right] \right).$$
\end{enumerate}
\end{cor}

\proof Part one is immediate from the Transfer Formula. Part two
results from the Shift Formula since the Roman shift of a
polynomial is multiplication by $x.\Box$

\subsection{Lagrange Inversion}
The explicit formulas\footnote{In particular, eq.~\ref{seams}.} of the
preceeding section lead to an elegant proof of the Lagrange Inversion
formula for the compositional inverse $f^{(-1)}(t) $ of a delta series
$f(t) $.

\begin{thm}[Lagrange Inversion] \label{Large} Let $f(t)$ be a delta
series, and $g(t)$ be any Laurent series of degree $d$, then
\begin{equation}\label{range}
g(f^{(-1)}(t)) = \sum_{k=d}^\infty
\aug{g(t)f'(t)f(t)^{-1-k}\frac{1}{x}}_{(1)} t^k.
\end{equation}
In particular,
$$ f^{(-1)}(t)^n = \sum_{k=n}^\infty
\aug{t^nf'(t)f(t)^{-1-k}\frac{1}{x}}_{(1)} t^k.$$
In other words, the coefficient of $t^k$ in $f^{(-1)}(t)^n$ is the
coefficient of $t^{n-1}$ in $f'(t)f(t)^{-1-k}$.
\end{thm}

\proof It will suffice to demonstrate eq.~\ref{range}. By eq.~\ref{seams}
formula,
$$ \r{k}! \aug{g(\D) p_k^{(1)} (x)}_{(1)}= \aug{g(\D) f'(\D) f(\D)^
{-1-k}\left(\frac{1}{x}\right)}_{(1)}$$
where $p_k^\alpha (x)$ is the logarithmic sequence of binomial type
associated with the delta operator $f(\D)$. Now,
$$\aug{f(\D)^n p_k^\alpha (x)}_\alpha = \delta_{nk} \r{n}! = \aug{\D^n
\lambda_k^\alpha (x)}_\alpha,$$
so in an augmentation of an operator acting on a $p_k^\alpha(x)$, one
can substitute $f^{(-1)}(\D)$ for $\D$ in the operator provided one
also substitutes $\lambda_k^\alpha(x)$ for $p_k^\alpha(x)$. In
particular,
\begin{equation}\label{before}
\r{k}! \aug{g(\D) p_k^{(1)} (x)}_{(1)} = \r{k}!
\aug{g(f^{(-1)}(\D))p_k^{(1)}(x)}_{(1)}.
\end{equation}
However, by the Expansion theorem (Theorem~\ref{let}),
$\r{k}!\aug{g(f^{(-1)}(\D))p_k^{(1)}(x)}_{(1)}$ is the coefficient of
$\D^k$ in $g(f^{(-1)}(\D))$.
Putting all this information together, we have eq.~\ref{range}.$\Box$

\section{Examples}
\subsection{The Logarithmic Lower Factorial Sequence}

We now return to our example of a logarithmic sequence of binomial
type (Section~\ref{eggs}). The logarithmic lower factorial
$(x)_n^\alpha$ is the logarithmic sequence of binomial type associated
with the delta operator $f(\D) = \Delta = e^{\D} - I$.

By Theorem~\ref{residual}, to determine the residual series, we must
first compute the Pincherle derivative
$$f'(\D) = \Delta ' = e^\D \ = E^1.$$
Thus, the residual series is
$$(x)_{-1}^{(1)} = E^{1}\frac{1}{x} = \frac{1}{x+1}.$$
This completes our calculations of section~\ref{eggs}, since all of
the other terms $(x)_n^\alpha$ can be computed via the formula
$$(x)_n^\alpha = \r{n}!^{-1} \mbox{skip}_{\alpha, (1)} \D^{-n-1}
\left( \frac{1}{x+1} \right). $$
Howver, we should note that in some cases the recursion formula is
much more useful
$$(x)_n^\alpha = \sigma (x-1)_{n-1}^\alpha.$$
In particular, for $\alpha=(0)$ and $n\geq 0$  or $\alpha=(1)$ and
$n<0$, this is equivalent to
$$(x)_n^\alpha = x(x-1)_{n-1}^\alpha.$$
In these cases, we immediately have
$$(x)_n = x(x-1)\cdots (x-n+1)$$
for $n$ positive, and
$$(x)_n^{(1)} = \frac{1}{(x+1)(x+2)\cdots (x-n)}$$
for $n$ negative.

\subsection{The Laguerre Logarithmic Sequence}\label{soon}
Our next example,
the Laguerre logarithmic sequence $L_n^\alpha (x)$ and the Laguerre
polynomials $L_n(x)= L_n^{(0)}(x)$, are the logarithmic sequence and
polynomials sequence of binomial type associated with the Laguerre
operator $f(\D) = W$ defined by eq.~\ref{seven} or eq.~\ref{september}.

We have for all $n$ and $\alpha$ by the Transfer formula,
\begin{eqnarray}
L_n^\alpha (x) &=& f'(\D) \left(\frac{\D}{f(\D)}\right)
\lambda_n^\alpha (x) \nonumber \\
&=& -(\D - I)^{n-1}\lambda_{n}^\alpha (x) \label{rod} \\
&=& \sum_{k=0}^\infty (-1)^{n+k} {n-1 \choose k}
\frac{\r{n}!}{\r{n-k}!} \lambda_{n-k}^\alpha (x).\nonumber
\end{eqnarray}
In particular, for $\alpha= (0)$ and $n$ positive, we have
$$L_n(x)= \sum_{k=0}^{n-1}\infty (-1)^{n+k} {n-1 \choose k}
\frac{{n}!}{{n-k}!} x^{n-k}.$$

Note that for the Laguerre logarithmic sequence, we not only have
$$\mbox{skip}_{(0),(1)}L_n^{(1)}(x)=L_n(x)$$
 as required by eq.~\ref{skipped}, but we also have
\begin{equation}\label{special}
\mbox{skip}_{(1),(0)}L_n(x)=L_n^{(1)}(x)
\end{equation}
for $n$ positive.\footnote{The logarithmic sequences of the form $a^n
\lambda_n^\alpha (x)$ (for example the harmonic logarithm) are even
more special since they obey eq.~\ref{special} even for $n=0$.
Nevertheless, there are no logarithmic sequences satisfying eq.
\ref{special} for any negative integers $n$.} This is a very special
property; it is only true of logarithmic sequences associated with
delta operators of the form $a\D / (\D - b)$ for $a,b$ real scalars
($a\neq 0$).

However,
$$ L_0^{(1)} = \log x + \frac{1}{x} - \frac{1}{x^2} + \frac{2}{x^3} -
\frac{6}{x^4} + \cdots  $$
which is not equal to $\mbox{skip} _{(1),(0)}L_0(x) =\mbox{skip} _{(1),(0)}1 =
\log x $.

From eq.~\ref{rod} and the classical identity (for polynomials)
\begin{equation}\label{roger}
e^x \D ^{n} e^{-x} = (\D - I)^n,
\end{equation}
we have the {\em Rodrigues formula} for the Laguerre sequence
$$L_n (x) = -e^x\D ^{n-1}e^{-x}x^n.$$

\subsection{The Abel Logarithmic Sequence}\label{liable}
We continue our series of examples with the logarithmic generalization
$A_n^\alpha (x;b)$ of the Abel polynomials $A_n(x;b)$ associated with
the delta operator $f(\D) = E^a \D $ mentioned in Section~\ref{able}.

By the shift formla, for $n$ neither 0 nor 1, we have
$$A_n^\alpha (x) = \sigma \lambda_{n-1}^\alpha  (x-nb).$$
For $\alpha = (0)$ and $n>1$ or $\alpha=(1)$ and $n<0$, this
immediately gives us
$$x(x-nb)^{n-1}.$$
In particular, the residual series is
$$A_{-1}^{(1)}(x) = x(x+b)^{-2}.$$

It remains now only to compute the terms of degree zero and one.
The series of degree zero is computed via the transfer formula, and
turns out be much simpler than might have been expected:
\begin{eqnarray*}
A_0^{(1)}(x) &=& (\D E^{b})' E^{-b} \log x\\
&=& (I + b\D) \log x\\
&=& \log x + b/x.
\end{eqnarray*}

As an application, we infer from the logarithmic binomial theorem (eq.
\ref{b})
$$A_0^\alpha (x+a) = \sum_{k=0}^\infty \rc{0}{k} A_k(a;b)
A_k^{(1)}(x)$$
the remarkable identity
$$\frac{b}{x+a} + \log (x+a) = \frac{b}{x} + \log x +
\sum_{k=1}^\infty \frac{(-1)^{k+1}a(a-bk)^{k-1}x}{k(x+bk)^{k+1}}.$$
For example, we can substitute $a=1$, $b=2$, and $x=5$, and compute
the first twelve terms of both sides yielding equality to seven
decimal places.

In general, the transfer formula gives us
\begin{eqnarray*}
A_n^{\alpha}(x) &=& E^{-nb} (I+b\D) \lambda_n^\alpha (x)\\
&=& \lambda_n^\alpha(x-nb) + b\r{n} \lambda_{n-1}^\alpha (x-nb).
\end{eqnarray*}
For example, the term of degree one is given by
$$A_1^{(1)}(x)= x\log(x-a) +a-x.$$

Finally, consider the Taylor's theorem for the logarithimc Abel
sequence (Theorem~\ref{tail}).
Any $p(x)\in {\cal L}^+$ can be expanded in terms of the logarithmic
Abel sequence as follows
$$p(x) = \sum_{\alpha \neq (0)} \sum_{n=-\infty}^{N_\alpha}
\frac{a_n^\alpha}{\r{n}!} A_n^\alpha (x)$$
where $a_n^\alpha = \aug{E^{nb} \D^n p(x)}_\alpha$. In particular, for
$p(x)=\log x$, $a_n^\alpha = 0$ for $\alpha \neq (1)$ or $n$ positive,
but for $n$ nonpositive and $\alpha = (1)$ we have
\begin{eqnarray*}
c_{n}^{\alpha } &=& \aug{E^{nb}D^n \log x}_{(1)}\\
&=& \aug{E^{nb}\lambda_{-n}^{(1)}(x)/(-n)!}_{(1)}\\
&=& \aug{E^{nb}\lambda_{-n}^{(0)}(x)/(-n)!}_{(0)}\\
&=& \aug{(x+nb)^{-n}}_{(0)}/(-n)!\\
&=& (nb)^{-n}/(-n)!.
\end{eqnarray*}
Thus,
\begin{eqnarray}
\log x &=& \log x + \frac{b}{x} + \sum_{n=1}^\infty \frac{(-nb)^n}{n!
\r{-n}!} A_{-n}^{(1)}(x)\nonumber \\
-\frac{b}{x} &=&  \sum_{n=1}^\infty \frac{(-nb)^n(-1)^{n-1}}{n}
x(x+nb)^{-n-2} \nonumber \\
\frac{b}{x^2} &=&  \sum_{n=1}^\infty \frac{(nb)^n}{n}
(x+nb)^{-n-2} .\label{wrong}
\end{eqnarray}
Note that eq.~\ref{wrong} is the correct version of the calculations
in \cite[p. 105]{LA}. We thank Richard Askey for bringing our error to
our attention.

\subsection{The Logarithmic Ramey Sequences}
Let $p_n^\alpha$ and $p_n(x)=p_n^{(0)}(x)$ be the logarithmic and
polynomial sequences of binomial type associated with a delta operator
$f(\D)$. The sequences $p_n^\alpha(x;b)$ and $p_n(x;b)$ associated
with the delta operator $E^b f(\D)$ are called the logarithmic and
polynomial {\em Ramey sequences}\/\footnote{Also called the
Abelization.} of $p_n(x)$ and $p_n^\alpha(x)$. See Ramey's paper on 
polynomial Ramey sequences for $f(\D)$. For example, if
$f(\D) = \D$, the Ramey sequences are the Abel sequences dealt with in
sections~\ref{able} and~\ref{liable}.

The following proposition is very useful in the calculation of Ramey
sequences.

\begin{prop} Let $p_n^\alpha(x)$ and $p_n(x)$ be the Ramey sequences
for the logarithmic and polynomial sequences of binomial type
$q_n^\alpha(x) $ and $q_n(x)=q_n^{(0)}(x)$ for the delta operator $f(\D)$.
\begin{enumerate}
\item Then for $n$ neither zero
nor one, the explicit
relationship between $p_n(x)$ and $q_n(x)$ is given by
$$p_n(x)= \frac{xq_n(x-nb)}{x-nb}.$$
\item For all $n$, we have
$$p_n^\alpha(x) = q_n^\alpha (x-nb) + \r{n}b f'(\D)^{-1}
q_{n-1}^\alpha (x-nb).$$
\item For $n$ negative and $\alpha = (1)$, we have
$$p_n^{(1)}(x) = q_n^{(1)}(x-nb) + \r{n} b x^{-1} q_n^{\alpha}(x).$$
\end{enumerate}
\end{prop}

\proof (1) Immediate from part two of Corollary~\ref{cute}.

(2) Immediate from part one of Corollary~\ref{cute}.

(3) Since the Roman shift is multiplication by $x$ in this context, we
have by the Recurrence Formula
$$ x^{-1}q_n^\alpha(x) = f'(\D) q_{n-1}^\alpha (x).$$
The conclusion now immediately follows from part two.$\Box$

For example, let us now compute the polynomial Ramey sequence for the lower
factorial. This is the  polynomial sequences
 $G_n(x)$ commonly called {\em the Gould polynomials} associated
with the delta operator $\D \Delta$.
We have
$$G_n(x) = \frac{x(x-nb)_n}{x-nb} = x(x-nb-1)_{n-1}.$$

\section{Connection Constants}
Let $p_n^\alpha(x)$ and $q_n^\alpha (x)$ be logarithmic sequences
with
\begin{equation}\label{cough}
q_n^\alpha (x) = \sum_{k=-\infty}^n c_{nk}^\alpha \lambda_k^\alpha
(x).
\end{equation}
Then the sequence
$$r_n^\alpha (x) = \sum_{k=-\infty}^n c_{nk}^\alpha p_k^\alpha
(x)$$
is called the {\em umbral composition} of the logarithmic sequences
$p_n^\alpha(x)$ and $q_n^\alpha(x)$. To denote umbral composition, we
write
$$r_n^\alpha (x) = q_n^\alpha({\bf p}).$$

Restricting out attention to polynomials, we see the if $p_n(x)$ and
$q_n(x)$ with
$$q_n(x) = \sum_{k=0}^n c_{nk}x^k,$$
then the umbral composition is uniquely defined, and is given by
$$q_n({\bf p}) = \sum_{k=0}^n c_{nk} p_n(x).$$

When $p_n^\alpha(x)$ and $q_n^\alpha (x)$ are of binomial type, then
umbral composition has a remarkable property.
\begin{thm}
Let $p_n^\alpha(x)$ and $q_n^\alpha(x)$ be the logarithmic sequences
of binomial type associated with the delta operators $f(\D)$ and
$g(\D)$ respectively. Then their umbral composition
$r_n^\alpha(x)=p_n^\alpha({\bf q})$ is the logarithmic sequence of
binomial type associated with the delta operator $f(g(\D))$.
\end{thm}

\proof We must show the $r_n^\alpha(x)$ obeys eqs.~\ref{one} and
\ref{two}.  Let the coefficients of $p_n^\alpha (x)$ and $f(\D)$ be
given by
\begin{eqnarray*}
p_n^\alpha(x) &=& \sum_{k=-\infty}^n c_{nk} \lambda_k^\alpha(x)\\
f(\D) &=& \sum_{j=1}^\infty a_j \D^j
\end{eqnarray*}
We write $c_{nk}$ since by eq.~\ref{skipped} the coefficients do not
depend on $\alpha$.

To prove eq.~\ref{one}, we proceed as follows:
\begin{eqnarray*}
\aug{r_n^\alpha ({\bf q})}_{\alpha} &=& \aug{\sum_{k=-\infty}^n c_{nk}
q_k^\alpha(x)}_\alpha \\
 &=& \sum_{k=-\infty}^n c_{nk} \aug{q_k^\alpha(x)}_\alpha \\
 &=& \sum_{k=-\infty}^n c_{nk} \delta_{0k} \\
 &=& c_{n0}\\
 &=& \aug{p_n^\alpha(x)}_\alpha \\
 &=& \delta_{n0}.
\end{eqnarray*}

Next, it remains to prove eq.~\ref{two}.
\begin{eqnarray*}
f(g(\D))r_n^\alpha (x) &=& \sum_{j=1}^\infty \sum_{k=-\infty}^n a_j
c_nk g(\D)^j q_k^\alpha (x)\\
&=& \sum_{j=1}^\infty \sum_{k=-\infty}^n a_j
c_nk \frac{\r{k}!}{\r{k-j}!} q_{k-j}^\alpha (x).
\end{eqnarray*}
However,
\begin{eqnarray*}
\r{n} p_{n-1}^\alpha (x) &=& f(\D) p_n^\alpha (x) \\
&=& \sum_{j=1}^\infty \sum_{k=-\infty}^n a_j
c_nk \D^j \lambda_k^\alpha (x)\\
&=& \sum_{j=1}^\infty \sum_{k=-\infty}^n a_j
c_nk \frac{\r{k}!}{\r{k-j}!} \lambda_{k-j}^\alpha (x),
\end{eqnarray*}
so that
$$f(g(\D)) r_n^\alpha (x) = \r{n} r_{n-1}^\alpha (x).\Box$$

By restricting our attention to ${\cal L}^{(0)}$, we have the
analogous result for polynomials.
\begin{cor}
Let $p_n(x)$ and $q_n(x)$ be the polynomial sequences
of binomial type associated with the delta operators $f(\D)$ and
$g(\D)$ respectively. Then their umbral composition
$r_n(x)=p_n({\bf q})$ is the  sequence of
binomial type associated with the delta operator $f(g(\D)).\Box$
\end{cor}

The preceeding theorem allows us to solve the {\em problem of
connection constants}. That is, given two logarithmic  sequence of
binomial type  $r_n^\alpha
(x)$ and $q_n^\alpha (x)$ associated with the delta operators $g(\D)$
and $h(\D)$, we are to find constants\footnote{Such
constants exists since the two sequences are basis, and the
coefficients do not depend on $\alpha$ by eq.~\ref{skipped}.} $c_{nk}$
such that
$$r_n^\alpha (x) = \sum_{-\infty}^n c_{nk} q_k^\alpha (x).$$
The preceeding theory shows that the constants $c_{nk}$ are given by a
logarithmic sequence $p_n^\alpha (x)$ as in eq.~\ref{cough} associated
with the delta operator $f(\D)=h(g^{(-1)}(\D))$.
Of course, the same result holds true for polynomial sequences of
binomial type.

In particular, given any logarithmic sequence of binomial type
$p_n^\alpha(x)$ associated with the delta operator $f(\D)$, there
exists a unique inverse logarithmic sequence
$\overline{p}_n^{\alpha}(x)$ associated with the delta operator
$f^{(-1)}(\D)$. In view of the logarithmic Taylor's theorem (Theorem
\ref{tailor}),
$$\overline{p}_n^\alpha (x) = \sum_{k=-=infty}^n \frac{\aug{\D^k
\overline{p}_n^\alpha(x)}}{\r{k}!} \lambda_k^\alpha(x).$$
However, by the demonstration of eq.~\ref{before},
$\overline{p}_n^\alpha(x)$ is the conjugate logarithmic sequence for the
delta operator $f(\D)$; that is,
$$\overline{p}_n^\alpha(x)=\sum_{k=-\infty}^n \frac{\aug{f(\D)^k
\lambda_n^\alpha(x)}_\alpha}{\r{n}!} \lambda_k^\alpha(x).$$

Moveover, by further umbral composition
$$\overline{p}_n^\alpha({\bf q}) = \sum_{k=-\infty}^n \frac{\aug{f(\D)^k
\lambda_n^\alpha(x)}_\alpha}{\r{n}!} q_k^\alpha(x).$$

Obviously, similar results hold for polynomials.

\section{Examples}
\subsection{From the Logarithmic Lower Factorial to the Harmonic
Logarithm}

To compute the connection constants expressing the harmonic logarithms
$\lambda_n^\alpha(x)$ in terms of the logarithmic lower factorials
$(x)_n^\alpha$ via the above techniques, we must first calculate the
logarithmic generalization $\phi_n^\alpha(x)$ of the {\em exponential
polynomials}  $\phi_n(x)$ associated with the delta operator
$\Delta^{(-1)}= \log (I+\D)$ and which have been studied by Touchard \cite{T}
among others.

As noted above, this is the logarithmic
conjugate sequence for forward difference operator $\Delta$. That is,
$$\phi_n^\alpha(x)=\sum_{k=-\infty}^n \frac{\aug{\Delta^k
\lambda_n^\alpha(x)}_\alpha}{\r{n}!} \lambda_k^\alpha(x).$$
In particular, for $\alpha = (0)$,
$$\phi_n(x)=\sum_{k=0}^n \frac{\aug{\Delta^k x^n}_\alpha}{\r{n}!}
x^k = \sum_{k=0}^n S(n,k) x^k$$
where $S(n,k)$ denote the {\em Stirling numbers of the second kind.}

To calculate the residual series, we as usual employ Theorem
\ref{residual},
\begin{eqnarray*}
\phi_{-1}^{(1)}(x) &=& \log(I+\D)' x^{-1}\\
&=& \sum_{k=0}^\infty (-1)^k \D^k x^{-1}\\
&=& \sum_{k=0}^\infty k!/x^{k+1}
\end{eqnarray*}
which is equivalent via umbral composition to eq.~\ref{golden}.

Other terms of the logarithmic exponential sequence are most easily
computed via the recurrence formula.
\begin{eqnarray*}
\phi_n^\alpha (x) &=& \sigma (\log(I+\D)')^{-1} \phi_{n-1}^\alpha(x)\\
&=& \sigma (I+\D) \phi_{n-1}^\alpha(x).
\end{eqnarray*}
For $\alpha=(0)$ and $n$ positive, or $\alpha=(1)$ and $n$ negative,
we thus have by eq.~\ref{roger}
\begin{eqnarray*}
\phi_n^\alpha(x) &=& xe^{-x}\D e^x \phi_{n-1}^\alpha(x)\\
&=& e^{-x}(x\D)e^x \phi_{n-1}^\alpha(x).
\end{eqnarray*}
In particular, for polynomials, we have by induction
$$\phi_n(x) = e^{-x}(x\D)^{n}e^x.$$

\subsection{From the Lower Factorials to the Upper Factorials}

The {\em logarithmic upper factorial} $\aug{x}_n^\alpha$ is the
logarithmic sequences associated with the backward difference
operator $\nabla = I-E^{-D}$. In effect, it is the important special
case of the Gould sequence with $b=-1$. Thus, for $\alpha=(0)$ and $n$
positive, we have the sequence of polynomials
$$\aug{x}_n = \aug{x}_n^\alpha = x(x+1)\cdots (x+n-1),$$
and for $\alpha=(1)$ and $n$ negative, we have
$$\aug{x}_n^{(1)} = \frac{1}{(x-1)(x-2)\cdots(x+n)}.$$

To express $(x)_n^\alpha$ in terms of $\aug{x}_n^\alpha$, we must
first compute the logarithmic sequence of binomial type
$p_n^\alpha(x)$ associated
with the delta operator $f(\D) = \Delta(\nabla^{(-1)})$.
\begin{eqnarray*}
f(\D) &=& \Delta(\nabla^{(-1)}) \\
 &=& \exp(-\log(I-\D) ) - I \\
 &=& \frac{I}{I-\D} - I \\
 &=& \frac{\D}{I-\D} \\
 &=& -W
\end{eqnarray*}
where $W$ is the Weierstrass operator. Hence,
the logarithmic sequence $p_n^\alpha(x)$ is essential the logarithmic
Laguerre sequence.
$$p_n^\alpha(x) = (-1)^n L_n^\alpha (x) = \sum_{k=0}^\infty (-1)^{k}
{n-1 \choose k} \frac{\r{n}!}{\r{n-k}!} \lambda_{n-k}^\alpha (x).$$
Thus,
$$(x)_n^\alpha= \sum_{k=0}^\infty (-1)^{k}
{n-1 \choose k} \frac{\r{n}!}{\r{n-k}!} \aug{x}_{n-k}^\alpha.$$
In particular, for $\alpha = (0)$ and $n$ positive,
$$x(x-1)\cdots(x-n+1) = \sum_{k=0}^{n-1} (-1)^{k}
{n-1 \choose k} \frac{{n}!}{{n-k}!} x(x+1)\cdots(x+n-k-1)$$
a fact which was not at all obvious {\em a priori}.

\subsection{From the Abel Series to the Harmonic Logarithms}

To compute the harmonic logarithms $\lambda_n^\alpha(x)$ in terms of
the logarithmic Abel sequence $A_n^\alpha(x;b)$ it will suffice to
compute the logarithmic inverse Abel sequence
$\mu_n^\alpha(x;b)=\overline{A}_n^\alpha(x;b)$ associated with the
composition inverse of $\D E^b$.

Unfortunately, there is no closed form
expression for the compositional inverse of $\D \exp(b\D)$.
Nevertheless, we can compute $\mu_n^\alpha(x)$ using the fact that it
is the conjugate logarithmic sequence for $\D E^b$.
$$\mu_n^\alpha(x)=\sum_{k=-\infty}^n \frac{\aug{\D^k E^{bk}
\lambda_n^\alpha(x)}_\alpha}{\r{n}!} \lambda_k^\alpha(x).$$
Now, $\aug{\D^k E^{bk} \lambda_n^\alpha(x)}_\alpha$ is equal to
${(bk)^{n-k}\r{n}!}/{\r{n-k}!}$.
Thus,
$$\mu_n^\alpha(x;b)=\sum_{k=-\infty}^n \frac{(bk)^{n-k}}{\r{n}!}
\lambda_k^\alpha(x).$$
And therefore
$$\lambda_n^\alpha(x)=\sum_{k=-\infty}^n \frac{(bk)^{n-k}}{\r{n}!}
A_k^\alpha(x;b),$$
or for $\alpha=(0)$ and $n$ positive or $\alpha=(1)$ and $n$ negative,
$$x^n = \sum_{k=-\infty}^n \frac{(bk)^{n-k}}{\r{n}!}
x(x-nb)^{n-1}.$$

\end{document}